\documentclass[12pt]{amsart}

\usepackage{amsthm}
\usepackage{amssymb}
\usepackage{amsmath}
\usepackage{mathtools}  
\numberwithin{equation}{section}
\usepackage{graphics,graphicx,colortbl}
\usepackage{subfigure}
\usepackage{multicol}
\usepackage{epsfig}
\usepackage{graphicx}
\usepackage[all,knot]{xy}
\xyoption{arc}
\usepackage{txfonts}
\definecolor{Mygrey}{gray}{0.8}
\newcommand{\bea}{\begin{eqnarray}}
\newcommand{\eea}{\end{eqnarray}}
\newcommand{\be}{\begin{eqnarray*}}
\newcommand{\ee}{\end{eqnarray*}}

\begin{document}
\title[Inverses of Cartan Matrices]{Inverses of Cartan Matrices of Lie Algebras and Lie Superalgebras}
\author[Yangjiang Wei]{Yangjiang Wei}
\address{School of Mathematics and Statistics, Guangxi Teachers Education University, Nanning 530023, P. R. China}
\email{gus02@163.com}
\author[Yi Ming Zou]{Yi Ming Zou}
\address{Department of Mathematical Sciences, University of Wisconsin, Milwaukee, WI 53201, USA} \email{ymzou@uwm.edu}
\thanks{* Corresponding author. Email: ymzou@uwm.edu}
\maketitle
\begin{abstract}
We derive explicit formulas for the inverses of the Cartan matrices of the simple Lie algebras and the basic classical Lie superalgebras, as well as for their infinite generalizations.
\end{abstract}
\section{Introduction}
\par
Cartan matrices show up in the context of Lie algebras, Lie groups, and quantum groups, they encode the defining relation information. In the finite dimensional case, the basic Cartan matrices are the ones correspond to the simple Lie algebras and the basic classical Lie superalgebras. Although one may consider the simple Lie algebras as part of the basic classical Lie superalgebras \cite{Kac1}, we separate them to make a distinction between the non-supercase and the supercase. The inverses of Cartan matrices also appear in a number of related topics \cite{Chin, Clin, Hue, KP, Malk, Mr, Qua, S1, SL}. However, except the treatment of the inverses of Cartan matrices in \cite{Lu}, there does not seem to be readily accessible explicit information on the inverses of Cartan matrices in spite of the fact that they can be described by quite nice formulas. For example, the Cartan matrix of type $A_n =(a_{ij})$ ($n\ge 1$):
\be
a_{ij} = \begin{cases} 2 &\mbox{if}\; i = j, \\
                                   -1 &\mbox{if}\; |i-j| = 1,\\
                                   0 &\mbox{otherwise},
            \end{cases}
\ee
has the following formula for the entries of its inverse:
\be
(A^{-1}_n)_{ij} = \min\{i,j\}-\frac{ij}{n+1}.
\ee 
\par
Although formulas for the inverses of the Cartan matrices of the finite dimensional simple Lie algebras like the one above can be obtained by using row operations and induction,  the method runs into some messy calculations in the supercase. Here, we use certain matrix update methods to derive explicit formulas for the inverses by starting with matrices whose inverses can be found easily, and thus avoiding messy calculations. This approach generalizes to other cases, allows us to derive the inverses of the Cartan matrices of the basic classical Lie superalgebras, as well as infinite generalizations of the Cartan matrices for both the Lie algebra case and the Lie superalgebra case, in a systematic way. The formulas for the finite cases are provided in Section 2, and the proofs for these formulas are given in Section 3. The infinite generalizations of the Cartan matrices and their inverses, both for the simple Lie algebras and the basic classical Lie superalgebras of types $A$, $B$, $C$, and $D$, are provided in Section 4. The exceptional cases are listed in the Appendix for completeness.
\section{Inverses of Cartan Matrices}
Since the inverses of the Cartan matrices for the exceptional Lie algebras $E_6$, $E_7$, $E_8$, $F_4$, $G_2$, and the exceptional Lie superalgebras $D(2,1;\alpha)$, $F(4)$, and $G(3)$  can be readily computed, we provide the formulas for the series $A$, $B$, $C$, and $D$ in this section, and list the exceptional cases in the Appendix.  The inverses of the Cartan matrices are given by explicit formulas for their entries in this section, the proofs in Section 3 contain the formulas in block matrix form.
\par\medskip
{\bf 2.1. Cartan Matrices of Lie Algebras $A, B, C, D$ and Their Inverses} 
\par\medskip
We use the Cartan matrix list in \cite{Hum}. Since $C_n = B_n^T$, we only list the inverse for $B_n$. The Cartan matrices and their inverses are:
\par\medskip
\be
 &A_n = \begin{pmatrix}
  2 & -1 & {} & {} & {} \\
  -1 & 2 & -1 & {} & {} \\
  {} & -1 & \ddots &\ddots &{} &{} \\
  {} & {} &\ddots & 2  & -1 \\
  {} & {} & {} & -1 & 2
\end{pmatrix} \; (n\ge 1), \\
&{}\\
& (A_n^{-1})_{ij} = \min\{i,j\} - \frac{ij}{n+1},\; 1\le i, j\le n.
\ee
\be
 B_n &=& \begin{pmatrix}
  2 & -1 & {} & {} & {} \\
  -1 & 2 & -1 & {} & {} \\
  {} & -1 & \ddots &\ddots &{} &{} \\
  {} & {} &\ddots & 2  & -2 \\
  {} & {} & {} & -1 & 2
\end{pmatrix} \; (n\ge 2), \\
{}&{}\\ 
(B_n^{-1})_{ij} &=& \frac{\min\{i,j\}}{1-\min\{0,n-i-1\} } =\begin{cases}
                                                                                     \min\{i,j\}, &\mbox{if}\; i < n, \\
                                                                                     \frac{j}{2}, &\mbox{if}\; i = n, \\
                                                                                    \end{cases} \; 1\le i, j\le n.
\ee
\be
 D_n &=& \begin{pmatrix}
  2 & -1 & {} & {} & {} \\
  -1 & \ddots & \ddots & {} & {} \\
  {} & \ddots & 2 & -1 &-1 \\
  {} & {} &-1& 2  & 0 \\
  {} & {} & -1 & 0 & 2
\end{pmatrix} \; (n\ge 4), \\
{}&{}\\
(D_n^{-1})_{ij} &=& \begin{cases} 
                                    i, &\mbox{if}\; 1\le i \le j \le n-2, \\
                                   i/2, &\mbox{if}\; i< n -1, j = n-1\;{\mbox{or}}\; n,\\
                                   \frac{n-2}{4}, &\mbox{if}\; i =n-1, j=n,\\
                                   n/4, &\mbox{if}\; i=j=n-1\;{\mbox{or}}\; n.
            \end{cases}
\ee
\par\medskip
Note that in the case of $D_n^{-1}$, we only listed the entries for $i\le j$ since $D_n^{-1}$ is symmetric. The proofs in Section 3 contain formulas in block matrix form for the inverses listed above. Note also that the formula holds for $n =2, 3$ -- this fact will be used for the supercase later, see the proof for the inverse of $D_{m,n}$.
\par\medskip
It can be seen from the formulas of the inverse matrices of $A_n, B_n$ ($C_n$), and $D_n$ that all entries of the inverses are positive. Since the entries of the inverses of the Cartan matrices for the exceptional simple Lie algebras are also positive, we conclude that the entries of the inverses of the Cartan matrices of simple Lie algebras are always positive. This fact was proved in \cite{Lu} using a different approach. We will see that the entries of the inverses of the Cartan matrices of Lie superalgebras can be positive, negative, and $0$.
\par\medskip
{\bf 2.2. Cartan Matrices of Lie Superalgebras $A(m,n)$, $B(m,n)$, $C(n)$, $D(m,n)$ and Their Inverses} 
\par\medskip
Unlike the Lie algebra case, Cartan matrices for the basic classical Lie superalgebras are not unique up to the equivalence defined by the Weyl group actions. For our purpose here, we will treat the Cartan matrices that correspond to the distinguished Borel subalgebras \cite{Kac2}. Using vectors $\epsilon_i$ and $\delta_j$ in a real vector space (the dimension is defined by the corresponding algebra)  with a non-degenerated bilinear form such that ($\delta_{ij} = 1$ if $i=j$ and $0$ otherwise)
\be
(\epsilon_i,\epsilon_j) = \delta_{ij},\; (\epsilon_i,\delta_j) = 0,\; (\delta_i,\delta_j) = -\delta_{ij},\;\forall\; i, j,
\ee
we express the distinguished simple root systems of $A(m,n)$, $B(m,n)$, $C(n)$,  and $D(m,n)$ as in \cite{Kac1, Kac2}  (we omit $A(n,n)$ since its Cartan matrix is singular):
\par
{\footnotesize
\be
{} & A(m,n)\; (m\ne n, m, n \ge 0): \{\epsilon_1 - \epsilon_2,  \ldots, \epsilon_{m+1}-\delta_1, \delta_1-\delta_2,\ldots, \delta_n-\delta_{n+1} \};\\
{}&{}\\
{}& B(m,n)\; (m, n>0): \{\delta_1-\delta_2, \ldots, \delta_n-\epsilon_1, \epsilon_1 - \epsilon_2, \ldots, \epsilon_{m-1}-\epsilon_m, \epsilon_{m} \};\\
{}&{}\\
{}&B(0,n)\; (n \ge 2): \{\delta_1-\delta_2,\delta_2-\delta_3, \ldots, \delta_n\};\\
{}&{}\\
{}&C(n)\; (n\ge 2): \{\epsilon_1 - \delta_1, \delta_1-\delta_2,\delta_2-\delta_3, \ldots, \delta_{n-2}-\delta_{n-1},2\delta_{n-1} \};\\
{}&{}\\
{}& D(m,n)\; (m\ge 2, n> 0): \{\delta_1-\delta_2,\ldots, \delta_n-\epsilon_1, \epsilon_1 - \epsilon_2, \ldots, \epsilon_{m-1}-\epsilon_m, \epsilon_{m-1}+\epsilon_m \}.
\ee}
\par
For $A(m,n)$, the distinguished Cartan matrix (of size $(m+n+1)\times (m+n+1)$) and its inverse are:
\be
  A_{m,n} =\left( \begin{array}{r|r|l}
  A_{m} & -1 & {}  \\
  \hline
  -1 & 0 & 1  \\
  \hline
   {} & 1& -A_n 
  \end{array}\right), \quad
  A_{m,n}^{-1} = \begin{pmatrix}
  L_1 & L_2  \\
  L_3 & L_4 
  \end{pmatrix}, \;\; m\ne n,\; m, n\ge 0,
\ee
where $A_m$ ($A_n$) is an $m\times m$  ($n\times n$) Cartan matrix of type $A$, $L_1$ is an $(m+1)\times (m+1)$ block, $L_2$ is an $(m+1)\times n$ block, $L_3 = L_2^T$, and $L_4$ is an $n\times n$ block, defined by
\be
(L_1)_{ij} &=& \min\{i,j\}+\frac{ij}{n-m},\; 1\le i,j\le m+1,\\
(L_2)_{ij} &=& \frac{i(n+1-j)}{n-m},\; 1\le i\le m+1, 1\le j\le n,\\
(L_4)_{ij} &=& \frac{(m+1)(n+1-i-j) +ij}{n-m} - \min\{i,j\},\; 1\le i,j\le n.
\ee
\par\medskip
For $B(m,n)$, since the distinguished Cartan matrix for $B(0,n)$ is $-B_n$, the negative of the Cartan matrix of the Lie algebra $B_n$, we will give the formula for the cases where $m, n >0$. The distinguished Cartan matrix (of size $(m+n)\times (m+n)$) and its inverse are:
\be
  B_{m,n} &=& \left( \begin{array}{r|l|l}
  -A_{n-1} & 1 & {}  \\
  \hline
  1 & 0 & -1  \\
  \hline
   {} & -1& B_m 
  \end{array}\right),\;\; m, n > 0, \\
  {}&{}\\
  (B_{m,n}^{-1})_{ij} &=& \begin{cases} 
                                    -\min\{i,j\}, &\mbox{if}\; 1\le i \le n, 1\le j \le m+n \\
                                    {}              &\mbox{or}\; 1\le j\le n, n+1\le i < m+n,\\
                                   \min\{i,j\}-2n, &\mbox{if}\; n<i< m+n, n<j \le m+n,\\
                                   -\frac{j}{2}, &\mbox{if}\; i =m+n, 1\le j\le n,\\
                                   \frac{j}{2}-n, &\mbox{if}\; i=m+ n, n< j\le m+n.    \end{cases}
\ee
\par\medskip
For $C(n)$, $n\ge 2$, we have:
\be
  C_{1,n-1} &=& \left( \begin{array}{r|l}
  0 & 1  \\
  \hline
  1 & -C_{n-1} 
  \end{array}\right) = -B_{n-1,1}^T,
  \ee
so
\be
(C_{1,n-1}^{-1})_{ij} = \begin{cases}
                                    2-\min\{i,j\}, &\mbox{if}\; 1\le j< n,\\
                                    1-\frac{i}{2}, &\mbox{if}\; j = n.
                                   \end{cases}
\ee
\par
For $D(m,n)$, we distinguish between the case $n=1$ and the case $n>1$ to make the formulas (and their proofs) more clear. We will list the entries such that $i\le j$ for the inverses since they are symmetric matrices. The distinguished Cartan matrix (of size $(m+1)\times (m+1)$) for $D(m,1)$ and its inverse are:
\par
\be
   D_{m,1} &=& \left( \begin{array}{r|l}
  0 & -1  \\
  \hline
  -1 & D_{m} 
  \end{array}\right),\;\; m\ge 2, \\
  {}&{}\\
  (D_{m,1}^{-1})_{ij} &=& \begin{cases} 
                                    i-2, &\mbox{if}\; 1\le i\le j \le m-1, \\
                                   \frac{i}{2}-1, &\mbox{if}\; 1\le i< m,\; j = m\; \mbox{or}\;m+1,\\
                                   \frac{m-3}{4}, &\mbox{if}\; i=m,  j = m+1,\\
                                   \frac{m-1}{4}, &\mbox{if}\; i =j=m\;\mbox{or}\; i=j=m+1.
                                       \end{cases}
\ee
\par
{\flushleft The distinguished Cartan matrix (of size $(m+n)\times (m+n), n>1$) for $D(m,n)$ and its inverse are:}
\par
\be
   D_{m,n} &=& \left( \begin{array}{r|l|l}
  -A_{n-1} & 1 &{}  \\
  \hline
  1 & 0 & -1 \\
  \hline
  {} & -1 & D_m
  \end{array}\right), \;\; m, n > 1,\\
  {}&{}\\
  (D_{m,n}^{-1})_{ij} &=& \begin{cases} 
                                    -i, &\mbox{if}\; 1\le i\le n,  i\le j \le m+n-2, \\
                                   -\frac{i}{2}, &\mbox{if}\; 1\le i \le n, j=m+n-1\;\mbox{or}\; m+n,\\
                                   i-2n,  &\mbox{if}\; n< i\le j < m+n-1,\\
                                   \frac{i}{2}-n,  &\mbox{if}\; n\le i < m+n-1, j=m+n-1\;\mbox{or}\; m+n,\\
                                   \frac{m-n-2}{4}, &\mbox{if}\; i =m+n-1, j=m+n,\\
                                   \frac{m-n}{4}, &\mbox{if}\; i=j = m+n-1\;\mbox{or}\; m+n.    \end{cases}
\ee
\par\medskip
\section{Proofs of the Formulas}
\par
We begin with the following auxiliary $n\times n$ symmetric matrix and its inverse:
\be
& S_n = \begin{pmatrix}
  2 & -1 & {} & {} & {} \\
  -1 & 2 & -1 & {} & {} \\
  {} & -1 & \ddots &\ddots &{} &{} \\
  {} & {} &\ddots & 2  & -1 \\
  {} & {} & {} & -1 & 1
\end{pmatrix} \; (n\ge 1), \\
&{}\\
 &(S_n^{-1})_{ij} = \min\{i,j\},\; 1\le i, j\le n.
\ee
\par
We can derive the inverse $S_n^{-1}$ by using elementary row operations on the matrix $(S_n|I_n)$ to reduce $S_n$ to the identity matrix. We first go upward by adding row $n$ to row $n-1$, then adding row $n-1$ to row $n-2$, and so on, to get
\be 
\left( \begin{array}{cccc|cccc}
  1 & {} & {} & {}  & 1 & 1 & \cdots & 1\\
  -1 & 1 & {} & {} & {} & 1 & \cdots & 1\\
  {} & \ddots & \ddots &{} & {} & {} &\ddots &\vdots\\
  {} & {} & -1 & 1  & {} & {} &{} & 1
\end{array}\right).
\ee
Then we go downward by adding row $1$ to row $2$ and so on to get
\par
\be 
\left( \begin{array}{ccccc|ccccc}
  1 & {} & {} & {}  &{} & 1 & 1 & 1 & \cdots & 1\\
  {} & 1 & {} & {} &{} & 1 & 2 & 2 & \cdots & 2\\
  {} &{} & 1 &{} &{} &   1 & 2 & 3 & \cdots & 3\\
  {} & {} &{} & \ddots &{} & \vdots & \vdots & \vdots &\ddots &\vdots\\
  {} & {} & {} & {}  & 1 & 1 & 2 & 3 & \cdots & n
\end{array}\right).
\ee
That gives the desired formula for $S_n^{-1}$.
\par
Now we can use the following rank-one updating formula \cite{Hag} (which can be verified directly) to find the inverses of $A_n$ and $B_n$: If $A$ is an $n\times n$ invertible matrix and $c, d$ are $n\times 1$ columns such that $1 + d^TA^{-1}c \ne 0$, then $A+cd^T$ is invertible and its inverse is given by
\bea\label{E1}
(A + cd^T)^{-1} = A^{-1} - \frac{A^{-1}cd^TA^{-1}}{1 + d^TA^{-1}c}.
\eea
\par
Let $E_{ij}$ be the matrix units (with $1$ at the $ij$-spot and $0$ elsewhere) and let $e_i$ be the the standard column vectors (with $1$ at the $i$-spot and $0$ elsewhere), $1\le i, j\le n$. We introduce the following column vectors:
\be
&(\underrightarrow{\mathbf{n}}) =(1,2,\ldots, n)^T, \quad (\underleftarrow{\mathbf{n}}) =(n,n-1,\ldots, 1)^T,\\
&(\mathbf{1_n}) =(\underbrace{1,1,\ldots, 1}_{\mbox{$n$ copies}})^T.
\ee
Note that $(\underrightarrow{\mathbf{n}})$ is the $n$-th column of $S_n^{-1}$.
\par\medskip
{\bf Remark.} In the proofs below, we frequently use the facts that $Ae_i$ is the $i$-column of the matrix $A$, $e_i^TA$ is the $i$-row of $A$, and $e_i^TAe_j$ is the $(i,j)$-entry of $A$.
\par\medskip
{\bf 3.1. Proofs for the Lie Algebra Cases} 
\par\medskip
{\bf Proof for $A_n^{-1}$}. Since
\be
A_n = S_n + E_{nn} = S_n + e_ne_n^T \quad\mbox{and}\quad 1+e_n^TS_n^{-1}e_n = 1+n,
\ee
by (\ref{E1}) we have :
\be
A_n^{-1} = S_n^{-1} - \frac{S_n^{-1}e_ne_n^TS_n^{-1}}{1+n} = S_n^{-1}  - \frac{1}{1+n}(\underrightarrow{\mathbf{n}})(\underrightarrow{\mathbf{n}})^T,
\ee
where the second equal sign holds since $S_n^{-1}e_n$ is the $n$-th column of $S_n^{-1}$ and $e_n^TS_n^{-1}$ is the $n$-th row.
This leads to the formula for $A_n^{-1}$ in Section 2. Q.E.D.
\par
\par\medskip
{\bf Proof for $B_n^{-1}$}.  We have:
\be
B_n &=& S_n + E_{nn} -E_{n,n-1} = S_n + (e_n-e_{n-1})e_n^T, \quad\mbox{and}\\
 {}&{}& 1+e_n^TS_n^{-1}(e_n-e_{n-1}) = 2.
\ee
Thus by (\ref{E1})
\be
B_n^{-1} = S_n^{-1}  - \frac{S_n^{-1}(e_n-e_{n-1})e_n^TS_n^{-1}}{2} =S_n^{-1} - \frac{1}{2}e_n(\underrightarrow{\mathbf{n}})^T. \quad \mbox{Q.E.D.}
\ee
\par\medskip
{\bf Proof for $D_n^{-1}$}.  There are different ways one can derive the formula for $D_n^{-1}$. For example, one can use the rank-one updating formula (\ref{E1}) twice to derive the inverse. But since we already have the inverse formula for the Cartan matrix of type $A$, we can use a formula for the inverse of a matrix in block form \cite{Hag}. Let
\be
M = \left( \begin{array}{cc}
  T & U\\
  V & W 
  \end{array}\right) 
  \ee
and suppose that both $T$ and $S =W - VT^{-1}U$ are invertible, then
\bea\label{E2}
M^{-1} = \begin{pmatrix}
                T^{-1} + T^{-1}US^{-1}VT^{-1} & -T^{-1}US^{-1}\\
                -S^{-1}VT^{-1} & S^{-1}
             \end{pmatrix}.
\eea
Write
\be
D_n = \left( \begin{array}{cc}
  A_{n-1} & -e_{n-2}\\
  -e_{n-2}^T & 2
  \end{array}\right),
\ee
where $e_{n-2}$ is of size $(n-1)\times 1$ with the nonzero entry at the $(n-2)$-spot, then the $S$ in (\ref{E2}) is of size $1\times 1$:
\be
S = 2 - e_{n-2}^TA_{n-1}^{-1}e_{n-2} =\frac{4}{n}.
\ee
Let 
\be
\mathbf{d} = (1,2,\ldots, n-2,\frac{n-2}{2})^T, 
\ee
then $-T^{-1}US^{-1} = \frac{1}{2}\mathbf{d}$, and 
\be
D_n^{-1} = \begin{pmatrix}
                A_{n-1}^{-1} + \frac{1}{n} \mathbf{d}\mathbf{d}^T & \frac{1}{2}\mathbf{d}\\
                {}&{}\\
                \frac{1}{2}\mathbf{d}^T & \frac{n}{4}
             \end{pmatrix}. \qquad \mbox{Q.E.D.}
\ee
\par\medskip
{\bf 3.2. Proofs for the Lie Superalgebra Cases}
\par\medskip
In the Lie superalgebra cases, we introduce the following auxiliary $(m+1)\times (m+1)$ symmetric matrix:
\be
 & R_{m+1} = \begin{pmatrix}
  2 & -1 & {} & {} & {} \\
  -1 & 2 & -1 & {} & {} \\
  {} & -1 & \ddots &\ddots &{} &{} \\
  {} & {} &\ddots & 2  & -1 \\
  {} & {} & {} & -1 & 0
\end{pmatrix} =S_{m+1}-e_{m+1}e_{m+1}^T\; (m\ge 1).
\ee
Its inverse can be obtained by using the updating formula (\ref{E1}) and the formula for the entries of $S_{m+1}^{-1}$. The entries of its inverse are:
\be 
(R_{m+1}^{-1})_{ij} = \min\{i,j\} -\frac{ij}{m},\; 1\le i, j\le m+1.
\ee
\par\medskip
For a matrix $M$, we let $r_i(M)$ (resp. $c_j(M)$) denote its $i$-th row (resp. $j$-th column).
\par\medskip
{\bf Proof for $A(m,n)$}. We have:
\be
  A_{m,n} =\left( \begin{array}{r|r|l}
  A_{m} & -1 & {}  \\
  \hline
  -1 & 0 & 1  \\
  \hline
   {} & 1& -A_n 
  \end{array}\right) = \left(\begin{array}{c|l}
  R_{m+1} & 1  \\
  \hline
  0 & -A_n
  \end{array} \right) +e_{m+2}e_{m+1}^T.
\ee
Let 
\be
U =  \left(\begin{array}{c|c}
  R_{m+1} & e_{m+1}e_1^T  \\
  \hline
  0 & -A_n
  \end{array} \right), 
\ee
where $e_{m+1}$ is $(m+1)\times 1$ and $e_1$ is $n\times 1$. Then
\bea\label{E31}
U^{-1} =  \left(\begin{array}{c|c}
  R_{m+1}^{-1} & R_{m+1}^{-1}e_{m+1}e_1^TA_n^{-1}  \\
  \hline
  0 & -A_n^{-1}
  \end{array} \right),
  \eea
where 
\bea\label{E32}\qquad
R_{m+1}^{-1}e_{m+1}e_1^TA_n^{-1} = c_{m+1}(R_{m+1}^{-1})r_1(A_n^{-1}) =-\frac{1}{m(n+1)}(\underrightarrow{\mathbf{m+1}})(\underleftarrow{\mathbf{n}})^T.
\eea
Since
\be
1+e_{m+1}^TU^{-1}e_{m+2} = 1 + (U^{-1})_{m+1,m+2} = \frac{m-n}{m(n+1)},
\ee
the rank-one updating formula (\ref{E1}) gives
\be
A_{m,n}^{-1} &=& U^{-1} -\frac{m(n+1)}{m-n}U^{-1}e_{m+2}e_{m+1}^TU^{-1}\\
{} & {} &\\
{} & =& U^{-1} + \frac{m(n+1)}{n-m}c_{m+2}(U^{-1})r_{m+1}(U^{-1}).
\ee
From (\ref{E31}) and (\ref{E32}), we have
\be
c_{m+2}(U^{-1}) = \begin{pmatrix}
     c_1(R_{m+1}^{-1}e_{m+1}e_1^TA_n^{-1})\\ 
     {}\\
     c_1(-A_n^{-1})  \end{pmatrix} = \begin{pmatrix}
         -\frac{n}{m(n+1)}(\underrightarrow{\mathbf{m+1}})\\
         {}\\
         -\frac{1}{n+1}(\underleftarrow{\mathbf{n}})
         \end{pmatrix},
\ee
and
\be
r_{m+1}(U^{-1}) &=& (r_{m+1}(R_{m+1}^{-1}),r_{m+1}(R_{m+1}^{-1}e_{m+1}e_1^TA_n^{-1}))\\
        {} & {} &\\
         {} &=& 
         (-\frac{1}{m}(\underrightarrow{\mathbf{m+1}})^T,
         -\frac{m+1}{m(n+1)}(\underleftarrow{\mathbf{n}})^T).
\ee
Thus if we write
\be
A_{m,n}^{-1} = \begin{pmatrix}
                         L_1 & L_2\\
                         L_3 & L_4
                        \end{pmatrix},
\ee
then $L_3 = L_2^T$ and (view $c_{m+2}(U^{-1})$ as a $2\times 1$ block matrix and $r_{m+1}(U^{-1})$ as $1\times 2$)
\be
L_1 &=& R_{m+1}^{-1}+\frac{n}{m(n-m)}(\underrightarrow{\mathbf{m+1}})(\underrightarrow{\mathbf{m+1}})^T,\\
L_2 &=& \frac{1}{n-m}(\underrightarrow{\mathbf{m+1}})(\underleftarrow{\mathbf{n}})^T,\\
L_4 &=& -A_n^{-1}+\frac{m+1}{(n-m)(n+1)}(\underleftarrow{\mathbf{n}})(\underleftarrow{\mathbf{n}})^T.\qquad \mbox{Q.E.D.}
\ee
\par\medskip
{\bf Proof for $B(m,n)$}. We have
\be
  B_{m,n} &=& \left( \begin{array}{r|l|l}
  -A_{n-1} & 1 & {}  \\
  \hline
  1 & 0 & -1  \\
  \hline
   {} & -1& B_m 
  \end{array}\right) = \left(\begin{array}{c|l}
  -R_{n} & -1  \\
  \hline
  0 & B_m
  \end{array} \right) -e_{n+1}e_{n}^T.
  \ee
Let 
\be
V = \left(\begin{array}{c|l}
  -R_{n} & -1  \\
  \hline
  0 & B_m
  \end{array} \right) = \left(\begin{array}{c|l}
  -R_{n} & e_ne_1^T  \\
  \hline
  0 & B_m
  \end{array} \right),
\ee
where $e_n$ is $n\times 1$ and $e_1$ is $m\times 1$. Then (compare with (\ref{E31}) and (\ref{E32})) 
\be
V^{-1} = \left(\begin{array}{c|l}
  -R_{n}^{-1} & -c_n(R_n^{-1})r_1(B_m^{-1})  \\
  \hline
  0 & B_m^{-1}
  \end{array} \right),
\ee
where
\be
-c_n(R_n^{-1})r_1(B_m^{-1}) = \frac{1}{n-1}(\underrightarrow{\mathbf{n}})(\mathbf{1_m})^T.
\ee
Let
\be
\mathbf{b_m} &=&(\underbrace{1,\ldots, 1}_{\mbox{{\tiny$m-1$ copies}}},1/2)^T,\;m\ge 1.
\ee
Note that $\mathbf{b_m}$ is the first column of $B_m^{-1}$. Then by the rank-one updating formula (\ref{E1}), we have
\be
B_{m,n}^{-1} = \left(\begin{array}{c|l}
  -R_{n}^{-1} -\frac{1}{n-1}(\underrightarrow{\mathbf{n}})(\underrightarrow{\mathbf{n}})^T & 
  -(\underrightarrow{\mathbf{n}})(\mathbf{1_m})^T  \\
  {}&{}\\
  \hline\\
  -\mathbf{b_m}(\underrightarrow{\mathbf{n}})^T & B_m^{-1} -n\mathbf{b_m}(\mathbf{1_m})^T
  \end{array} \right).
\ee
Note that for the above computations to make sense, we need $n>1$, but the final formula for the inverse does not involve dividing by $n-1$ and holds for all $n>0$ (see the formula for the inverse of $B_{m,n}$ in Section 2). Q.E.D.
\par\medskip
{\bf Proof for $D(m,n)$}. For the case $n =1$, we have
\be
   D_{m,1} = \left( \begin{array}{r|l}
  0 & -1  \\
  \hline
  -1 & D_{m} 
  \end{array}\right) = D_{m+1} - 2e_1e_1^T,
\ee
so the rank-one updating formula (\ref{E1}) gives
\be
D_{m,1}^{-1} = D_{m+1}^{-1} - 2\mathbf{d_{m+1}}\mathbf{d_{m+1}}^T,
\ee
where
\be
\mathbf{d_{m+1}} &=&(\underbrace{1,\ldots, 1}_{\mbox{{\tiny$m-1$ copies}}},1/2,1/2)^T, \;m\ge 2,
\ee
is the first column of $D_{m+1}^{-1}$.
\par
If $n>1$, then
\be
  D_{m,n} &=& \left(  \begin{array}{r|l|l}
  -A_{n-1} & 1 &{}  \\
  \hline
  1 & 0 & -1 \\
  \hline
  {} & -1 & D_m
  \end{array}\right) = \left(\begin{array}{c|l}
  -R_{n} & -1  \\
  \hline
  0 & D_m
  \end{array} \right) -e_{n+1}e_{n}^T,
\ee
thus an argument similar to that of the case $B_{m,n}$ gives
\be
D_{m,n}^{-1} = \left(\begin{array}{c|l}
  -R_{n}^{-1} -\frac{1}{n-1}(\underrightarrow{\mathbf{n}})(\underrightarrow{\mathbf{n}})^T & 
  -(\underrightarrow{\mathbf{n}})\mathbf{d_m}^T  \\
  {}&{}\\
  \hline\\
  -\mathbf{d_m}(\underrightarrow{\mathbf{n}})^T & D_m^{-1} -n\mathbf{d_m}\mathbf{d_m}^T
  \end{array} \right).    
\ee
From this block matrix and the fact that the formula for the inverse of $D_m$ holds for $m \ge 2$, we obtain the entries of $D_{m,n}^{-1}$. Q.E.D.
\par
\section{Infinite Cases}
The Cartan matrices of types $A, B, C$ and $D$ can be generalized to infinite case naturally (see pp.112-114 in \cite{Kac3} and references \cite{Le, S2}). We give the inverse formulas for these matrices and the corresponding supercases. One can derive these formulas from their finite counter parts by using ``limit'' process and taking into consideration of the re-labelling of the vertices in the diagrams (which corresponds to row and column indexing). For example, in the case of $A_{\infty}$, the inverse can be obtained by letting $n \rightarrow \infty$ in the formula for the inverse of $A_n$. For another example, we note that one can obtain the formula of the inverse of $A_{\infty,\infty}$ from that of $A_{m,\infty}$ by considering the fact that the special vertex labelled by $m+1$ in the diagram of $A_{m,\infty}$ is labelled as vertex $0$ in the diagram of $A_{\infty,\infty}$, so we set $m+1$ to $0$. The shifts in the formulas are due to the extra label $0$ in the diagram of $A_{\infty,\infty}$, because the label starts from $1$ for $A_{m,\infty}$. We omit the details since the formulas for the inverses can be verified directly. 
\par\medskip
For $A_{\infty}$, if
 \vspace{0.5cm}
   \begin{center}\setlength{\unitlength}{1cm}
    \begin{picture}(0,0)
    \put(-5,0){$A_{\infty}$:}
    \put(-3,0){$\medcirc$}
    \put(-2,0){$\medcirc$}
    \put(-1,0){$\medcirc$}
    \put(-2.75,0.1){\line(1,0){0.75}}
    \put(-1.75,0.1){\line(1,0){0.75}}
    \put(-0.75,0.1){\line(1,0){0.75}}
    \put(0.1,0.1){$\ldots$}
    \put(0.1,-0.4){$\ldots$}
    \put(-2.98,-0.5){{\footnotesize 1}}
    \put(-1.98,-0.5){{\footnotesize 2}}
    \put(-0.98,-0.5){{\footnotesize $3$}}
    \put(0.85,-0.45){{\footnotesize $\longrightarrow\infty$}}
    \put(1.9,-0.1){,}
    \put(2.2,0){then}
    \put(3, 0){ $(A_{\infty}^{-1})_{ij}=\min\{i,j\}$;}
    \end{picture}
   \end{center}\vspace{0.6cm}
and if
   \vspace{0.5cm}
   \begin{center}\setlength{\unitlength}{1cm}
    \begin{picture}(0,0)
    \put(-5.25,0){$A_{\infty}$:}
    \put(1,0){$\medcirc$}
    \put(0,0){$\medcirc$}
    \put(-1,0){$\medcirc$}
    \put(0.25,0.1){\line(1,0){0.75}}
    \put(-0.75,0.1){\line(1,0){0.75}}
    \put(-1.75,0.1){\line(1,0){0.75}}
    \put(-2.3,0.1){$\ldots$}
    \put(-2.3,-0.4){$\ldots$}
    \put(0.95,-0.5){{\footnotesize $-1$}}
    \put(-0.05,-0.5){{\footnotesize $-2$}}
    \put(-1.05,-0.5){{\footnotesize $-3$}}
    \put(-3.75,-0.45){{\footnotesize $-\infty\longleftarrow$}}
    \put(1.6,-0.1){,}
    \put(1.9,0){then}
    \put(2.8, 0){$(A_{\infty}^{-1})_{ij}=\min\{-i,-j\}$.}
    \end{picture}
   \end{center}\vspace{1cm}
   \par\medskip
   For
     \vspace{0.5cm}
   \begin{center}\setlength{\unitlength}{1cm}
    \begin{picture}(0,0)
    \put(-5.75,0){$B_{\infty}$:}
    \put(0.1,0){$\medcirc$}
    \put(-0.75,0){$\medcirc$}
    \put(-1.75,0){$\medcirc$}
    \put(-0.5,0){$\Longrightarrow$}
    \put(-1.5,0.1){\line(1,0){0.75}}
    \put(-2.5,0.1){\line(1,0){0.75}}
    \put(-3.15,0.1){$\ldots$}
    \put(-3.15,-0.4){$\ldots$}
    \put(0.15,-0.5){{\footnotesize $0$}}
    \put(-0.95,-0.5){{\footnotesize $-1$}}
    \put(-1.82,-0.5){{\footnotesize $-2$}}
    \put(-4.45,-0.45){{\footnotesize $-\infty\longleftarrow$}}
    \put(0.5,-0.1){,}
    \put(1, 0){$(B_{\infty}^{-1})_{ij}= \begin{cases} \frac{j}{2}, & i=0,\\
                                                                                              \min\{i,j\}, & i\ne 0. \end{cases}$}
    \end{picture}
   \end{center}\vspace{1cm}
   \par\medskip
   For
   \vspace{1cm}
   \begin{center}\setlength{\unitlength}{1cm}
    \begin{picture}(0,0)
    \put(-5,0){$D_{\infty}$:}
    \put(0.95,-0.7){$\medcirc$}
    \put(0.25,0){$\medcirc$}
    \put(-0.75,0){$\medcirc$}
    \put(0.95,0.7){$\medcirc$}
    \put(0.45,0.2){\line(1,1){0.5}}
    \put(0.5,0){\line(1,-1){0.5}}
    \put(-0.5,0.1){\line(1,0){0.75}}
    \put(-1.5,0.1){\line(1,0){0.75}}
    \put(-2.05,0.1){$\ldots$}
    \put(-2.05,-0.4){$\ldots$}
     \put(1.45,0.7){{\footnotesize $1$}}
    \put(1.45,-0.7){{\footnotesize $0$}}
    \put(0.1,-0.5){{\footnotesize $-1$}}
    \put(-0.8,-0.5){{\footnotesize $-2$}}
    \put(-3.5,-0.45){{\footnotesize $-\infty\longleftarrow$}}
    \put(2,0){,}
    \end{picture}
   \end{center}\vspace{1cm}
   and for $j\le i$ (the matrix is symmetric),
   \be
   (D_{\infty}^{-1})_{ij} = \begin{cases} \frac{1}{4}, & i= j = 0\;\mbox{or}\; 1,\\
                                                                      - \frac{1}{4}, & i=1, j = 0,\\    
                                                                       \frac{j}{2}, & j<0, i = 1\;\mbox{or}\; 0,\\  
                                                                         j, & j\le i\le -1. 
                                             \end{cases}
   \ee
   \par
   For
    \vspace{0.5cm}
      \begin{center}\setlength{\unitlength}{1cm}
       \begin{picture}(0,0)
       \put(-5.5,0){$A_{m,\infty}$:}
       \put(-4,0){$\medcirc$}
        \put(-3.93,-0.5){{\footnotesize $1$}}
       \put(-3.1,0){$\cdots$}
       \put(-2,0){$\medcirc$}
       \put(-1,0){$\otimes$}
       \put(-2.45,0.1){\line(1,0){0.45}}
       \put(-1.75,0.1){\line(1,0){0.75}}
       \put(-0.75,0.1){\line(1,0){0.75}}
        \put(0,0){$\medcirc$}
       \put(0.3,0.1){\line(1,0){0.75}}
        \put(-3.75,0.1){\line(1,0){0.45}}
        \put(1.25,0){$\cdots$}
       \put(1.25,-0.5){$\cdots$}
       \put(-3.1,-0.5){$\cdots$}
       \put(-1.98,-0.5){{\footnotesize $m$}}
       \put(-1.25,-0.5){{\footnotesize $m+1$}}
       \put(-0.2,-0.5){{\footnotesize $m+2$}}
       \put(2.25,-0.45){{\footnotesize $\longrightarrow\infty$}}
       \put(3.5,-0.1){,}
       \end{picture}
      \end{center}\vspace{0.6cm}
 the Cartan matrix and its inverse are:
 \be
   A_{m,\infty} =\left( \begin{array}{r|r|l}
   A_{m} & -1 & {}  \\
   \hline
   -1 & 0 & 1  \\
   \hline
    {} & 1& -A_{\infty} 
   \end{array}\right), \quad
   A_{m,\infty}^{-1} = \begin{pmatrix}
   L_1 & L_2  \\
   L_3 & L_4 
   \end{pmatrix}, 
 \ee
 where $L_1$ is an $(m+1)\times (m+1)$ block, $L_2$ is an $(m+1)\times \infty$ block, $L_3 = L_2^T$, and $L_4$ is an $\infty\times\infty$ block, defined by
 \be
 (L_1)_{ij} &=& \min\{i,j\},\; 1\le i,j\le m+1,\\
 (L_2)_{ij} &=& i,\; 1\le i\le m+1, 1\le j,\\
 (L_4)_{ij} &=& m+1 - \min\{i,j\},\; 1\le i,j.
 \ee
 \par\medskip
    For
     \vspace{0.5cm}
       \begin{center}\setlength{\unitlength}{1cm}
        \begin{picture}(0,0)
        \put(-5.5,0){$A_{\infty,\infty}$:}
        \put(-2.5,0){$\cdots$}
        \put(-2.5,-0.5){$\cdots$}
        \put(1.03,0){$\medcirc$}
        \put(0,0){$\otimes$}
        \put(-1.75,0.1){\line(1,0){0.75}}
        \put(-0.75,0.1){\line(1,0){0.75}}
        \put(0.25,0.1){\line(1,0){0.75}}
         \put(-1,0){$\medcirc$}
        \put(1.3,0.1){\line(1,0){0.75}}
         \put(-4,-0.5){$-\infty\longleftarrow$}
         \put(2.25,0){$\cdots$}
        \put(2.25,-0.5){$\cdots$}
        \put(-1.07,-0.5){{\footnotesize $-1$}}
        \put(0.03,-0.5){{\footnotesize $0$}}
        \put(1.1,-0.5){{\footnotesize $1$}}
        \put(3,-0.45){{\footnotesize $\longrightarrow\infty$}}
        \put(4,-0.1){,}
        \end{picture}
       \end{center}\vspace{0.6cm}
  the Cartan matrix and its inverse are:
  \be
    A_{\infty,\infty} =\left( \begin{array}{r|r|l}
    A_{\infty} & -1 & {}  \\
    \hline
    -1 & 0 & 1  \\
    \hline
     {} & 1& -A_{\infty} 
    \end{array}\right), \quad
    A_{\infty,\infty}^{-1} = \begin{pmatrix}
    L_1 & L_2  \\
    L_3 & L_4 
    \end{pmatrix}, 
  \ee
where $L_1$, $L_2$, and $L_4$ are $\infty\times \infty$ blocks, and $L_3 = L_2^T$, defined by
  \be
  (L_1)_{ij} &=& \min\{i-1,j-1\},\; i,j\le 0,\\
  (L_2)_{ij} &=& i - 1,\; i\le 0, 1\le j,\\
  (L_4)_{ij} &=&  - \min\{i+1,j+1\},\; 1\le i,j.
  \ee
  \par\medskip
  For
       \vspace{0.5cm}
         \begin{center}\setlength{\unitlength}{1cm}
          \begin{picture}(0,0)
          \put(-5.5,0){$B_{m,\infty}$:}
          \put(-2.5,0){$\cdots$}
          \put(-2.5,-0.5){$\cdots$}
          \put(1.03,0){$\medcirc$}
          \put(0,0){$\otimes$}
          \put(-1.75,0.1){\line(1,0){0.75}}
          \put(-0.75,0.1){\line(1,0){0.75}}
          \put(0.25,0.1){\line(1,0){0.75}}
           \put(-1,0){$\medcirc$}
          \put(1.35,0.1){\line(1,0){0.45}}
           \put(-4,-0.5){$-\infty\longleftarrow$}
           \put(2.05,0){$\cdots$}
          \put(2.75,0.1){\line(1,0){0.45}}
          \put(-1.07,-0.5){{\footnotesize $-1$}}
          \put(0.03,-0.5){{\footnotesize $0$}}
          \put(1.1,-0.5){{\footnotesize $1$}}
          \put(3.25,0){$\medcirc$}
           \put(3,-0.5){{\footnotesize $m-1$}}
            \put(3.55,0){$\Longrightarrow$}
            \put(4.25,0){$\medcirc$}
            \put(4.25,-0.5){{\footnotesize $m$}}
          \put(5,-0.1){,}
          \end{picture}
         \end{center}\vspace{0.6cm}
 \be
   B_{m,\infty} &=& \left( \begin{array}{r|l|l}
   -A_{\infty} & 1 & {}  \\
   \hline
   1 & 0 & -1  \\
   \hline
    {} & -1& B_m 
   \end{array}\right),\;\; m > 0, \\
   {}&{}\\
   (B_{m,\infty}^{-1})_{ij} &=& \begin{cases} 
                                     -\min\{i,j\}, &\mbox{if}\; i \le 0\; \mbox{or}\; j \le 0, \\
                                    \min\{i,j\}, &\mbox{if}\; 0<i< m, 0<j \le m,\\
                                    \frac{|j|}{2}, &\mbox{if}\; i=m.    \end{cases}
 \ee
 \par\medskip
   For
         \vspace{0.5cm}
           \begin{center}\setlength{\unitlength}{1cm}
            \begin{picture}(0,0)
            \put(-5.5,0){$D_{m,\infty}$:}
            \put(-2.5,0){$\cdots$}
            \put(-2.5,-0.5){$\cdots$}
            \put(1.03,0){$\medcirc$}
            \put(0,0){$\otimes$}
            \put(-1.75,0.1){\line(1,0){0.75}}
            \put(-0.75,0.1){\line(1,0){0.75}}
            \put(0.25,0.1){\line(1,0){0.75}}
             \put(-1,0){$\medcirc$}
            \put(1.35,0.1){\line(1,0){0.45}}
             \put(-4,-0.5){$-\infty\longleftarrow$}
             \put(2.05,0){$\cdots$}
            \put(2.75,0.1){\line(1,0){0.45}}
            \put(-1.07,-0.5){{\footnotesize $-1$}}
            \put(0.03,-0.5){{\footnotesize $0$}}
            \put(1.1,-0.5){{\footnotesize $1$}}
            \put(3.25,0){$\medcirc$}
             \put(2.75,-0.5){{\footnotesize $m-2$}}
              \put(4.15,0.7){$\medcirc$}
              \put(4.65,0.7){{\footnotesize $m$}}
               \put(4.15,-0.7){$\medcirc$}
               \put(3.55,0.2){\line(1,1){0.55}}
                \put(3.55,0){\line(1,-1){0.55}}
                 \put(4.65,-0.7){{\footnotesize $m-1$}}
            \put(5.75,0){,}
            \end{picture}
           \end{center}\vspace{0.6cm}
   \be
      D_{m,\infty} &=& \left( \begin{array}{r|l|l}
     -A_{\infty} & 1 &{}  \\
     \hline
     1 & 0 & -1 \\
     \hline
     {} & -1 & D_m
     \end{array}\right), \;\; m > 1,
     \ee
 and for $i\le j$ (the matrix is symmetric),  
   \be
     (D_{m,\infty}^{-1})_{ij} &=& \begin{cases} 
                                       -i, &\mbox{if}\; i\le 0,  i\le j \le m-2, \\
                                      -\frac{i}{2}, &\mbox{if}\; i \le 0, j=m-1\;\mbox{or}\; m,\\
                                      i,  &\mbox{if}\; 0< i\le j \le m-1,\\
                                      \frac{i}{2},  &\mbox{if}\; 0< i < m-1, j=m-1\;\mbox{or}\; m,\\
                                      \frac{m-2}{4}, &\mbox{if}\; i =m-1, j=m,\\
                                      \frac{m}{4}, &\mbox{if}\; i=j = m-1\;\mbox{or}\; m.    \end{cases}
   \ee
\par\medskip
\subsection*{Acknowledgment}
YJW is supported by the National Natural Science Foundation of China (11461010, 11661013, 11661014), the Guangxi Science Research and Technology Development Project (1599005-2-13), and  China Scholarship Council. YMZ acknowledges the support of a Simons Foundation Collaboration Grant for Mathematicians (416937).
 
\section*{Appendix}
We list the inverses of the Cartan matrices for the exceptional Lie algebras and Lie superalgebras in this appendix.
{\tiny
\be
E_6^{-1} = \begin{pmatrix}
                  2 & 0 & -1 & 0 & 0 & 0\\
                  0 & 2 & 0 & -1 & 0 & 0\\
                  -1 & 0 & 2 & -1 & 0 & 0\\
                   0 & -1 & -1 & 2 & -1 & 0\\
                   0 & 0 & 0 & -1 & 2 & -1\\
                   0 & 0 & 0 & 0 & -1 & 2
                  \end{pmatrix}^{-1}
                =  \frac{1}{3}\begin{pmatrix}
                  4 & 3 & 5 & 6 & 4 & 2\\
                  3 & 6 & 6 & 9 & 6 & 3\\
                  5 & 6 & 10 & 12 & 8 & 4\\
                  6 & 9 & 12 & 18 & 12 & 6\\
                  4 & 6 & 8 & 12 & 10 & 5\\
                  2 & 3 & 4 & 6 & 5 & 4
                  \end{pmatrix},
\ee
}
{\tiny
\be
E_7^{-1} = \begin{pmatrix}
                  2 & 0 & -1 & 0 & 0 & 0 & 0\\
                  0 & 2 & 0 & -1 & 0 & 0 & 0\\
                  -1 & 0 & 2 & -1 & 0 & 0 &0\\
                   0 & -1 & -1 & 2 & -1 & 0 & 0\\
                   0 & 0 & 0 & -1 & 2 & -1 & 0\\
                   0 & 0 & 0 & 0 & -1 & 2 &-1\\
                   0 & 0 & 0 & 0 & 0 & -1 & 2
                  \end{pmatrix}^{-1}
                =  \frac{1}{2}\begin{pmatrix}
                  4 & 4 & 6 & 8 & 6 & 4 & 2\\
                  4 & 7 & 8 & 12 & 9& 6 & 3\\
                  6 & 8 & 12 & 16 & 12 & 8 & 4\\
                  8 & 12 & 16 & 24 &18 & 12 & 6\\
                  6 & 9 & 12 & 18 & 15 & 10 & 5\\
                  4 & 6 & 8 & 12 & 10 & 8 & 4\\
                  2 & 3 & 4 & 6 & 5 & 4 & 3
                  \end{pmatrix},
\ee
}
{\tiny
\be
E_8^{-1} = \begin{pmatrix}
                  2 & 0 & -1 & 0 & 0 & 0 & 0 & 0\\
                  0 & 2 & 0 & -1 & 0 & 0 & 0 & 0\\
                  -1 & 0 & 2 & -1 & 0 & 0 &0 & 0 \\
                   0 & -1 & -1 & 2 & -1 & 0 & 0 & 0\\
                   0 & 0 & 0 & -1 & 2 & -1 & 0 & 0\\
                   0 & 0 & 0 & 0 & -1 & 2 &-1 & 0\\
                   0 & 0 & 0 & 0 & 0 & -1 & 2 & -1\\
                    0 & 0 & 0 & 0 & 0 & 0 & -1 & 2
                  \end{pmatrix}^{-1}
                =  \begin{pmatrix}
                  4 & 5 & 7 & 10 & 8 & 6 & 4 & 2\\
                  5 & 8 & 10 & 15 & 12 & 9& 6 & 3\\
                  7 & 10 & 14 & 20 & 16 & 12 & 8 & 4\\
                  10 & 15 & 20 & 30& 24& 18 & 12& 6\\
                  8 & 12 & 16 & 24 & 20 &15 & 10 & 5\\
                  6 & 9 & 12 & 18 & 15 & 12 & 8 & 4\\
                  4 & 6 & 8 & 12 & 10 & 8 & 6 & 3\\
                  2 & 3 & 4 & 6 & 5 & 4 & 3 & 2
                  \end{pmatrix},
\ee
}
{\tiny
\be
F_4^{-1} = \begin{pmatrix}
                  2 &  -1 & 0 & 0 \\
                  -1 & 2 & -2 & 0\\
                   0 & -1 & 2 &  -1\\
                   0 & 0 & -1 & 2
                  \end{pmatrix}^{-1}
                = \begin{pmatrix}
                  2 & 3 & 4 & 2\\
                  3 & 6 & 8 & 4\\
                 2 & 4 & 6 & 3\\
                 1 & 2 & 3 & 2
                  \end{pmatrix},
\ee
}
{\tiny
\be
G_2^{-1} = \begin{pmatrix}
                  2 &  -1 \\
                  -3 & 2
                  \end{pmatrix}^{-1}
                =  \begin{pmatrix}
                  2 &  1 \\
                  3 & 2
                  \end{pmatrix},
\ee
}
{\tiny
\be
D(2,1;\alpha)^{-1} = \begin{pmatrix}
                  0 &  1 &\alpha \\
                  -1 & 2 & 0\\
                  -1 & 0 & 2
                  \end{pmatrix}^{-1}
                =   \frac{1}{1+\alpha}\begin{pmatrix}
                 2 &  -1 & -\alpha \\
                  1 & \frac{\alpha}{2} & -\frac{\alpha}{2}\\
                  1 & -\frac{\alpha}{2} & \frac{\alpha}{2}
                  \end{pmatrix}, \; \alpha\ne -1, 0,
\ee
}
{\tiny
\be
SF_{4}^{-1} = \begin{pmatrix}
                  0 &  1 & 0 & 0 \\
                  -1 & 2 & -2 & 0\\
                   0 & -1 & 2 &  -1\\
                   0 & 0 & -1 & 2
                  \end{pmatrix}^{-1}
                = \frac{1}{3}\begin{pmatrix}
                  2 & -3 & -4 & -2\\
                  3 & 0 & 0 & 0\\
                 2 & 0 & 2 & 1\\
                 1 & 0 & 1 & 2
                  \end{pmatrix},
\ee
}
{\tiny
\be
G_3^{-1} = \begin{pmatrix}
                  0 &  1 & 0 \\
                  -1 & 2 & -3\\
                  0 & -1 & 2
                  \end{pmatrix}^{-1}
                = \frac{1}{2} \begin{pmatrix}
                  1 & -2 &  3 \\
                  2 & 0  & 0\\
                  1 & 0 & 1
                  \end{pmatrix}.
\ee
}
\end{document}